\documentclass[english]{paper}
\usepackage[T1]{fontenc}
\usepackage[latin9]{inputenc}
\usepackage{varioref}
\usepackage{float}
\usepackage{fancybox}
\usepackage{calc}
\usepackage{textcomp}
\usepackage{amsthm}
\usepackage{amsmath}
\usepackage{amssymb}
\usepackage{setspace}
\PassOptionsToPackage{normalem}{ulem}
\usepackage{ulem}

\makeatletter

\providecommand{\tabularnewline}{\\}

\numberwithin{equation}{section}
\numberwithin{figure}{section}

\makeatother

\usepackage{babel}
\begin{document}

\title{{\huge Hexile Sieve Analysis of Prime and Composite Integers.}}

\author{© Roger M. Creft.}

\maketitle
keywords : prime sieving, prime-counting function, diophantine equation,
combinatorics.
\begin{abstract}
Here we demonstrate a sieve for analysing primes and their composites,
using equivalence classes based on the modulo 6 return value as applied
to the Natural numbers $\mathbb{N}$. Five features of this 'Hexile'
sieve are reviewed. The first aspect, is that it narrows the search
for primes to one-third of $\mathbb{N}$. The second feature is that
we can obtain from the equivalence class formulae, a property of its
diophantine equations to distinguish between primes and composites
resulting from multiplication of these primes. Thirdly we can from
these diophantine formulations ascribe a non-random occurence to not
only the composites in the two equivalence classes but by default
and as a consequence : non-randomness of occurence to the resident
primes. Fourthly we develop a theoretical basis for sieving primes.
Of final mention is that the diophantine equations allows another
route to a prime counting function using combinatorics or numerical
analysis. 
\end{abstract}

\section{Introduction to Hexile sieving of $\mathbb{N}$ using modulo 6.\label{sec:Hexile-Sieving-of-Natural Numbers}}

It is not difficult to arrange $\mathbb{N}$ the set of Natural numbers
into six equivalence classes, on the basis of their modulo 6 return
value as shown overleaf in table \ref{tab:Illustrating-Hexile-Sieving}.

\begin{table}[H]
\begin{tabular}{|c|c|c|c|c|c|c|}
\hline 
\textbf{Hexile integer-level} : n  & $\mathbb{H}^{\mathrm{1}}$ & $\mathbb{H}^{\mathrm{2}}$ & $\mathbb{H}^{3}$ & $\mathbb{H}^{\mathrm{4}}$ & $\mathbb{H}^{\mathrm{5}}$ & $\mathbb{H}^{\mathrm{0}}$\tabularnewline
\hline 
\hline 
0 & 1 & 2 & 3 & 4 & 5 & 6\tabularnewline
\hline 
1 & 7 & 8 & 9 & 10 & 11 & 12\tabularnewline
\hline 
2 & 13 & 14 & 15 & 16 & 17 & 18\tabularnewline
\hline 
3 & 19 & 20 & 21 & 22 & 23 & 24\tabularnewline
\hline 
4 & 25 & 26 & 27 & 28 & 29 & 30\tabularnewline
\hline 
5 & 31 & 32 & 33 & 34 & 35 & 36\tabularnewline
\hline 
6 & 37 & 38 & 39 & 40 & 41 & 42\tabularnewline
\hline 
7 & 43 & 44 & 45 & 46 & 47 & 48\tabularnewline
\hline 
8 & 49 & 50 & 51 & 52 & 53 & 54\tabularnewline
\hline 
9 & 55 & 56 & 57 & 58 & 59 & 60\tabularnewline
\hline 
$\vdots$ & $\vdots$ & $\vdots$ & $\vdots$ & $\vdots$ & $\vdots$ & $\vdots$\tabularnewline
\hline 
\end{tabular}

\caption{\label{tab:Illustrating-Hexile-Sieving}Hexile equivalence class sieving
of the first 60 integers in $\mathbb{Z^{\mathrm{+}}}$}
\end{table}

Where the six equivalence classes are referenced as Hexile integer
classes symbolized as $\mathbb{H}^{\mathrm{k}}$\textbf{, }where in
set-builder notation we define the six classes thus for $\mathrm{k\in\left\{ 0,1,2,3,4,5\right\} }$
:

\paragraph*{
\begin{equation}
\mathbb{H^{\mathrm{0}}=\left[\mathrm{0}\right]=\left\{ \mathrm{6,12,18,...\infty}\right\} \equiv\left\{ \mathit{x|x=\left(6\times n\right)+0,n\in\mathbb{Z^{\mathrm{+}}\mathit{,n\geq0}}}\right\} }\label{eq:H0}
\end{equation}
\begin{equation}
\mathbb{H^{\mathrm{1}}=\left[\mathrm{1}\right]=\left\{ \mathrm{7,13,19,...\infty}\right\} \equiv\left\{ \mathit{x|x=\left(6\times n\right)+1,n\in\mathbb{Z^{\mathrm{+}}\mathit{,n\geq0}}}\right\} }\label{eq:H11}
\end{equation}
\begin{equation}
\mathbb{H^{\mathrm{2}}=\left[\mathrm{2}\right]=\left\{ \mathrm{8,14,20,...\infty}\right\} \equiv\left\{ \mathit{x|x=\left(6\times n\right)+2,n\in\mathbb{Z^{\mathrm{+}}\mathit{,n\geq0}}}\right\} }\label{eq:H2}
\end{equation}
\begin{equation}
\mathbb{H^{\mathrm{3}}=\left[\mathrm{3}\right]=\left\{ \mathrm{9,15,21,...\infty}\right\} \equiv\left\{ \mathit{x|x=\left(6\times n\right)+3,n\in\mathbb{Z^{\mathrm{+}}\mathit{,n\geq0}}}\right\} }\label{eq:H3}
\end{equation}
\begin{equation}
\mathbb{H^{\mathrm{4}}=\left[\mathrm{4}\right]=\left\{ \mathrm{10,16,22,...\infty}\right\} \equiv\left\{ \mathit{x|x=\left(6\times n\right)+4,n\in\mathbb{Z^{\mathrm{+}}\mathit{,n\geq0}}}\right\} }\label{eq:H4}
\end{equation}
\begin{equation}
\mathbb{H^{\mathrm{5}}=\left[\mathrm{5}\right]=\left\{ \mathrm{11,17,23,...\infty}\right\} \equiv\left\{ \mathit{x|x=\left(6\times n\right)+5,n\in\mathbb{Z^{\mathrm{+}}\mathit{,n\geq0}}}\right\} }\label{eq:H5}
\end{equation}
}

\bigskip{}

Closer examination of table \ref{tab:Illustrating-Hexile-Sieving},
reveals that $\mathbb{H}^{\mathrm{0}}$, $\mathbb{H}^{\mathrm{2}}$
and $\mathbb{H}^{\mathrm{4}}$ contains all the even numbers, consequently
with the exception of 2 in $\mathbb{H}^{\mathrm{2}}$, they are all
composites which are multiples of 2. We will also observe that $\mathbb{H}^{\mathrm{3}}$
contains only one element which is prime: 3, the other elements being
multiples of 3. 

Which leaves us with $\mathbb{H}^{\mathrm{1}}$ and $\mathbb{H}^{\mathrm{5}}$
as candidates for the domiciled sub-set partitioning of primes in
$\mathbb{N}$. Within which we will find all primes greater than 3,
as well as composites which are themselves product of these primes.

Hence prime numbers greater than or equal to 5, can be found in either
of the two HICs : 

$\mathbb{H^{\mathrm{1}}}$ and $\mathbb{H^{\mathrm{5}}}$ defined
above by equations \ref{eq:H11} and \ref{eq:H5}. %
\footnote{HIC is an abbreviation for Hexile integer class(es). %
}

\pagebreak{}

\section{Formulae of composites in $\mathbb{H^{\mathrm{1}}}$ and $\mathbb{H^{\mathrm{5}}}$.}

We will now examine the nature of composites in \emph{$\mathbb{H^{\mathrm{1}}}$
and $\mathbb{H^{\mathrm{5}}}$, }which also details the autochtnous
nature of their multiplicative origin from endemic primes. 

\rule[0.5ex]{1\columnwidth}{1pt}

\subsection{Composite resulting from multiplying two primes from \emph{$\mathbb{H^{\mathrm{1}}}$.\label{sub:Composite-resulting-from-c11}}}

\begin{singlespace}
If we were to consider two primes \emph{$\mathbb{\mathit{p_{1},q_{1}\:\in\:}H^{\mathrm{1}}}$}
with formulation \emph{$\mathbb{\mathit{\: p_{1}\:=\:\left(6\times m\right)\:+1}}$
}and\emph{ $\mathbb{\mathit{q_{1}\:=\:\left(6\times n\right)\:+1}}$}
where $\mathit{m\mathit{,n\:\in\:}\mathbb{Z^{\mathrm{+}}}}$.
\end{singlespace}

Multiplication of primes \emph{$\mathbb{\mathit{p_{1}}}$ }and \emph{$\mathbb{\mathit{q_{1}}}$}
yields a composite $c^{11}$\emph{ }whence :\emph{ 
\[
\mathbb{\mathit{\mathbb{\mathit{c^{11}\:=\: p_{1}\times q_{1}}}\mathbb{\mathit{\:\qquad\qquad\qquad\qquad\qquad\qquad}}}}
\]
\[
\mathbb{\mathit{=\:\mathbb{\mathit{\left[\left(6\times m\right)\:+1\right]}}\times\mathit{\left[\left(6\times n\right)\:+1\right]\:}}}
\]
\begin{equation}
\mathbb{\mathit{\;\quad\therefore\, c^{11}\:=\left(\mathbb{\mathit{36\times}\mathit{m\times n}}\right)}}\:+\:\mathit{\left(\mathit{6}\times m\right)\,+}\:\mathit{\left(\mathit{6}\times n\right)\:+\:1}\:\label{eq:Cpq11}
\end{equation}
}

Applying modulo 6 to the result of equation \ref{eq:Cpq11} returns
a value of 1, it then logically implies that the composite : $c^{11}$
is an element of \emph{$\mathbb{H^{\mathrm{1}}}$.}

\rule[0.5ex]{1\columnwidth}{1pt}

\subsection{Composite resulting from multiplying two primes from \emph{$\mathbb{H^{\mathrm{5}}}$.\label{sub:Composite-resulting-from-c55}}}

\begin{singlespace}
If we were to consider two primes \emph{$\mathbb{\mathit{p_{5},q_{5}\:\in\:}H^{\mathrm{5}}}$}
with formulations \emph{$\mathbb{\mathit{\: p_{5}\:=\:\left(6\times m\right)\:+5}}$}
and\emph{ $\mathbb{\mathit{q_{5}\:=\:\left(6\times n\right)\:+5}}$}
where $\mathit{m\mathit{,n\:\in\:}\mathbb{Z^{\mathrm{+}}}}$.
\end{singlespace}

Multiplication of primes \emph{$\mathbb{\mathit{p_{5}}}$ }and \emph{$\mathbb{\mathit{q_{5}}}$
}yields a composite $c^{55}$ whence :\emph{ 
\[
\mathbb{\mathit{\mathbb{\mathit{c^{55}\:=\: p_{5}\times q_{5}}}\mathbb{\mathit{\:\qquad\qquad\qquad\qquad\qquad\qquad}}}}
\]
}

\begin{singlespace}
\emph{
\[
\mathbb{\mathit{=\:\mathbb{\mathit{\left[\left(6\times m\right)\:+5\right]}}\times\mathit{\left[\left(6\times n\right)\:+5\right]\:}}}
\]
}

\emph{
\begin{equation}
\mathbb{\mathit{\;\quad\therefore\, c^{55}\:=\left(\mathbb{\mathit{36\times}\mathit{m\times n}}\right)}}\:+\:\mathit{\left(\mathit{30}\times m\right)\,+}\:\mathit{\left(\mathit{30}\times n\right)\:+\:25}\:\label{eq:Cpq11-1}
\end{equation}
}
\end{singlespace}

Applying modulo 6 to the result of equation \ref{eq:Cpq11-1} returns
a value of 1, it then logically implies that the composite : $c^{55}$
is an element of \emph{$\mathbb{H^{\mathrm{1}}}$.}

\rule[0.5ex]{1\columnwidth}{1pt}

\subsection{Composite resulting from multiplying one prime from \emph{$\mathbb{H^{\mathrm{1}}}$
}with another prime from \emph{$\mathbb{H^{\mathrm{5}}}$.\label{sub:Composite-resulting-from-c15}}}

\begin{singlespace}
If we were to consider two primes \emph{$\mathbb{\mathit{p_{1}\:\in\:}H^{\mathrm{1}}}$}
with formulation \emph{$\mathbb{\mathit{\: p_{1}\:=\:\left(6\times m\right)\:+1}}$}
and \emph{$\mathbb{\mathit{q_{5}\:\in\:}H^{\mathrm{5}}}$ }with formulation
\emph{$\mathbb{\mathit{q_{5}\:=\:\left(6\times n\right)\:+5}}$} where
$\mathit{m\mathit{,n\:\in\:}\mathbb{Z^{\mathrm{+}}}}$.
\end{singlespace}

Multiplication of primes \emph{$\mathbb{\mathit{p_{1}}}$ }and \emph{$\mathbb{\mathit{q_{5}}}$}
yields a composite $c^{15}$\emph{ }whence :\emph{ 
\[
\mathbb{\mathit{\mathbb{\mathit{\;\quad c^{15}\:=\: p_{1}\times q_{5}}}\mathbb{\mathit{\,=\mathbb{\mathit{\:\mathbb{\mathit{\left[\left(6\times m\right)\:+1\right]}}\times\mathit{\left[\left(6\times n\right)\:+5\right]}}}}}}}
\]
}

\begin{singlespace}
\emph{
\begin{equation}
\mathbb{\mathit{\;\quad\therefore\, c^{15}\:=\left(\mathbb{\mathit{36\times}\mathit{m\times n}}\right)}}\:+\:\mathit{\left(\mathit{30}\times m\right)\,+}\:\mathit{\left(\mathit{5}\times n\right)\:+\:5}\:\label{eq:Cpq11-1-1}
\end{equation}
}
\end{singlespace}

Applying modulo 6 to the result of equation \ref{eq:Cpq11-1-1} returns
a value of 5, it then logically implies that the composite : $c^{15}$
is an element of \emph{$\mathbb{H^{\mathrm{5}}}$.}

\rule[0.5ex]{1\columnwidth}{1pt}

\subsection{Possibility of \emph{$\mathbb{\mathit{\mathbb{H^{\mathrm{1}}}}}$
}and \emph{$\mathbb{\mathit{\mathbb{H^{\mathrm{5}}}}}$}composites
arising from multiplication of other HIC\emph{ }elements\emph{.}}

Elements from \emph{$\mathbb{\mathit{\;\mathbb{H^{\mathrm{0}}}}}$,
$\mathbb{\mathit{\mathbb{H^{\mathrm{2}}}}}$} and \emph{$\mathbb{\mathit{\;\mathbb{H^{\mathrm{4}}}}}$}
are all even numbers, since multiplying an even number with an odd
number always gives an even number, then their multiplication cannot
possibly give rise to \emph{$\mathbb{\mathit{\mathbb{H^{\mathrm{1}}}}}$
}and \emph{$\mathbb{\mathit{\mathbb{H^{\mathrm{5}}}}}$} composites.

We are left then with the sole possibility of an element of \emph{$\mathbb{\mathit{\mathbb{H^{\mathrm{3}}}}}$,}
giving rise to a composite in \emph{$\mathbb{\mathit{\mathbb{H^{\mathrm{1}}}}}$
}and \emph{$\mathbb{\mathit{\mathbb{H^{\mathrm{5}}}}}$.} Then for
\emph{$\mathbb{\mathit{r_{3}\in\mathbb{H^{\mathrm{3}}}}}$, }with
a generic formulation\emph{ }of\emph{ :}

\emph{
\begin{equation}
r_{3}=\left(\mathit{6\times}d\right)+3\label{eq:r-H3}
\end{equation}
}

We next need to examine the end-result of multiplying $r_{3}$ with
an element of \emph{$\mathbb{\mathit{\mathbb{H^{\mathrm{1}}}}}$,
$\mathbb{\mathit{\mathbb{H^{\mathrm{3}}}}}$} and \emph{$\mathbb{\mathit{\;\mathbb{H^{\mathrm{5}}}}}$.}
Table \ref{tab:Result-of-multiplying-H2} below, shows the results
of these three combinations, where \emph{$\mathbb{\mathit{p_{1}\in\mathbb{H^{\mathrm{1}}}}}$
}and \emph{$\mathbb{\mathit{q_{5}\:\in\:}H^{\mathrm{5}}}$} are formulated
thus\emph{: }

\emph{
\begin{equation}
\mathbb{\mathit{\: p_{1}\:=\:\left(6\times m\right)\:+1}}\label{eq:p1}
\end{equation}
\begin{equation}
\mathbb{\mathit{q_{5}\:=\:\left(6\times n\right)\:+}}\mathrm{5}\label{eq:q5}
\end{equation}
}

\begin{table}[H]
\begin{tabular}{|c|c|c|}
\hline 
 & $r_{3}=\left(\mathit{6\times}d\right)+3$ & modulo 6 of composite \tabularnewline
\hline 
\hline 
$r_{3}$ & $\mathbb{\mathit{\, c^{33}\:=\left(\mathbb{\mathit{36\times}\mathit{d}^{\mathrm{2}}}\right)}}\:+\:\mathit{\left(36\times d\right)\,+}\:9$ & 3\tabularnewline
\hline 
\emph{$\mathbb{\mathit{\: p_{1}\:}}$}  & $\mathbb{\mathit{\, c^{13}\:=\left(\mathbb{\mathit{36\times}\mathit{d\times m}}\right)}}\:+\:\mathit{\left(18\times m\right)\,+}\:\mathit{\left(6\times d\right)\:+3}\:$ & 3\tabularnewline
\hline 
\emph{$\mathbb{\mathit{q_{5}\:}}$}  & $\mathbb{\mathit{\, c^{35}\:=\left(\mathbb{\mathit{36\times}\mathit{d\times n}}\right)}}\:+\:\mathit{\left(\mathit{30}\times d\right)\,+}\:\mathit{\left(\mathit{18}\times n\right)\:+\:15}\:$ & 3\tabularnewline
\hline 
\end{tabular}\caption{\label{tab:Result-of-multiplying-H2}Results of multiplying \emph{$\mathbb{\mathit{r_{3}\in\mathbb{H^{\mathrm{3}}}}}$}
with an element of \emph{$\mathbb{\mathit{\;\mathbb{H^{\mathrm{1}}}}\mathrm{,}}$
}$\mathbb{\mathit{\mathbb{H^{\mathrm{3}}}}}$ and \emph{$\mathbb{\mathit{\;\mathbb{H^{\mathrm{5}}}}}.$}}

\end{table}

The results from table \ref{tab:Result-of-multiplying-H2}, indicates
an element of \emph{$\mathbb{\mathit{\mathbb{H^{\mathrm{3}}}}}$}
will not result in a\emph{ $\mathbb{\mathit{\mathbb{H^{\mathrm{1}}}}}$}
or \emph{$\mathbb{\mathit{\;\mathbb{H^{\mathrm{5}}}}}$ }composite
through multiplication. So then\emph{ }composites in \emph{$\mathbb{\mathit{\mathbb{H^{\mathrm{1}}}}}$}
and \emph{$\mathbb{\mathit{\;\mathbb{H^{\mathrm{5}}}}}$} are found
autochthonous with primes\emph{ }from which they are the product through
multiplication\emph{.}

\rule[0.5ex]{1\columnwidth}{1pt}

\section{Investigating the formulation of composites, through the contribution
of the Hexile Levels of their constituent primes.\label{sec:Investigating-the-formulation}}

Our probings will now take us to the inner sanctum of the formulaic
expression of composites, into - for want of a more apt term, we will
refer to as - the \emph{nucleus} of an element. 

\bigskip{}

\framebox{\begin{minipage}[t]{1\columnwidth}%

\subsubsection{Definition of the \emph{nucleus} value of a Natural number.\label{sub:Definition-of-the-Nucleus-Value}}

We define the \emph{nucleus} value of any Natural number, as the quotient
obtained when that number is divided by 6. %
\end{minipage}}

\medskip{}

Of note is that the \emph{nucleus }value can also be referred to as
the Hexile level of a Natural number, which was initially referenced
as the row numbers in the Hexile Sieve table \ref{tab:Illustrating-Hexile-Sieving}
of section (\ref{sec:Hexile-Sieving-of-Natural Numbers}).

Applying the above definition (\ref{sub:Definition-of-the-Nucleus-Value})
to the three formulations of \emph{$\mathbb{\mathit{\mathbb{H^{\mathrm{1}}}}}$
}and \emph{$\mathbb{\mathit{\mathbb{H^{\mathrm{5}}}}}$} composites
obtained in the previous sub-sections (\ref{sub:Composite-resulting-from-c11}),
(\ref{sub:Composite-resulting-from-c55}) and (\ref{sub:Composite-resulting-from-c15}),
we obtain the following three formulations for \emph{nuclei} values
as shown in table \vref{tab:Nucleus-value-Formualtion-H1-H5}.

\begin{table}[H]
\begin{tabular}{|c|c|}
\hline 
Composite Formulae & Nucleus Formulae\tabularnewline
\hline 
\hline 
$\mathbb{\mathit{c^{11}\:=\left(\mathbb{\mathit{36\times}\mathit{m\times n}}\right)}}+\left(\mathit{6}\times m\right)\,+\:\mathit{\left(\mathit{6}\times n\right)+1}\:$ & $\left(6\times m\times n\right)+\mathit{m+n}$\tabularnewline
\hline 
$\mathbb{\mathit{c^{55}\:=\left(\mathbb{\mathit{36\times}\mathit{m\times n}}\right)}}+\mathit{\left(\mathit{30}\times m\right)+}\mathit{\left(\mathit{30}\times n\right)+25}$ & $\left(6\times m\times n\right)+\left(5\times m\right)+\left(5\times n\right)+4$\tabularnewline
\hline 
$\mathbb{\mathit{c^{15}\:=\left(\mathbb{\mathit{36\times}\mathit{m\times n}}\right)}}+\left(\mathit{30}\times m\right)+\left(\mathit{5}\times n\right)+5$ & $\left(6\times m\times n\right)+\mathit{\left(5\times m\right)+\: n}$\tabularnewline
\hline 
\end{tabular}\caption{Showing nucleus value formulations for \emph{$\mathbb{\mathit{\mathbb{H^{\mathrm{1}}}}}$
}and \emph{$\mathbb{\mathit{\mathbb{H^{\mathrm{5}}}}}$} composites.\label{tab:Nucleus-value-Formualtion-H1-H5} }
\end{table}

\medskip{}

The results of table \vref{tab:Nucleus-value-Formualtion-H1-H5} can
be used as a preliminary investigative formulative-decomposition tool,
for an element of \emph{$\mathbb{\mathit{\mathbb{H^{\mathrm{1}}}}}$
}and \emph{$\mathbb{\mathit{\mathbb{H^{\mathrm{5}}}}}.$} 

Whereby if we are given a positive integer $c$, if its modulo 6 value
returned is 1 and its easily obtained \emph{nucleus} value is $Q$,
then we can ascribe to it as a composite one of the following two
formulations of equations \vref{eq:COLIDE-11} and \vref{eq:COLIDE-55}
given below :

\begin{equation}
\mathbb{\mathit{Q\,=\,}\mathit{\left(6\times m\times n\right)+\:\mathit{m+\:\mathit{n\:}\:\quad}\qquad}}\qquad\qquad\label{eq:COLIDE-11}
\end{equation}

or 

\begin{equation}
Q\,=\left(6\times m\times n\right)+\:\mathit{\left(5\times m\right)+\:\mathit{\left(5\times n\right)\,+\,4}}\label{eq:COLIDE-55}
\end{equation}

Also for a given Natural number $c$, if its modulo 6 value returned
is 5 and its \emph{nucleus} value is $Q$, we can specifically ascribe
to it as a composite the formulation of equation \vref{eq:COLIDE-15}
given below:
\begin{equation}
\mathit{\mathbb{\mathit{Q\,=\,}\mathit{\left(6\times m\times n\right)+\mathit{\left(5\times m\right)+\: n}}}\quad}\qquad\qquad\label{eq:COLIDE-15}
\end{equation}

Where $\mathit{m\mathit{,n\:\in\:}\mathbb{Z^{\mathrm{+}}}}$ are the
Hexile levels of the two constituent integers presumed to be primes
in the composite elements of \emph{$\mathbb{\mathit{\mathbb{H^{\mathrm{1}}}}}$
}and \emph{$\mathbb{\mathit{\mathbb{H^{\mathrm{5}}}}}$} as symbolized
by the variable $c$. 

\rule[0.5ex]{1\columnwidth}{1pt}

\section{Distinguishing feature of Composite linear diophantine equations
for elements of \emph{$\mathbb{\mathit{\mathbb{H^{\mathrm{1}}}}}$
}and \emph{$\mathbb{\mathit{\mathbb{H^{\mathrm{5}}}}}.$\label{sec:Distinguishing-feature-of-Composites_Primes}}}

We investigated in the previous section (\ref{sec:Investigating-the-formulation})
three diophantine formulations for composites in \emph{$\mathbb{\mathit{\mathbb{H^{\mathrm{1}}}}}$
}and \emph{$\mathbb{\mathit{\mathbb{H^{\mathrm{5}}}}}.$ }

These three equations viz : (\ref{eq:COLIDE-11}), (\ref{eq:COLIDE-55})
and (\ref{eq:COLIDE-15}), are essentially formulae defining the relationship
of a composite's \emph{nucleus} value in terms of the Hexile levels
of its two constituent primes/composites. 

\medskip{}

We will for ease of referential paltabilty designate equations (\ref{eq:COLIDE-11}),
(\ref{eq:COLIDE-55}) and (\ref{eq:COLIDE-15}) as \emph{\uline{Composite
Linear Diophantine Equation}}(s), or \emph{colide(s)} for short.

Specifically we will reference a\emph{ }colide-11\emph{ }to equation
\ref{eq:COLIDE-15} viz. :

\[
Q\,=\mathbb{\mathit{\,}\mathit{\left(6\times m\times n\right)}}\:+\:\mathit{\mathit{m}+\:\mathit{n\:\:}\:\;\quad}\qquad\qquad\ref{eq:COLIDE-15}
\]

where $\; m,\, n\:$ are the \textbf{HI-level}s of the two constituent
primes/composites with $\; m,\, n\:\in\:\mathbb{Z^{\mathrm{+}}},$
$m\,,n\:>\:0$, with $Q$ being the readily computable \emph{nucleus\index{Nucleus}}
value.\emph{ }

\medskip{}

Also a\emph{ }colide-55, we will reference\emph{ }to equation \ref{eq:COLIDE-55}
below :

\[
\mathit{\mathbb{\mathit{Q\,=\,}\mathit{\left(6\times m\times n\right)+\:\mathit{\left(5\times m\right)+\:\mathit{\left(5\times n\right)\,+\,4\:}\:\;\quad}}}\quad}\qquad\qquad\ref{eq:COLIDE-55}
\]

where $\; m,\, n$ are the \textbf{HI-level}s of the two constituent
\index{Primes}primes/composites\index{Composites} with $m,\, n\:\in\:\mathbb{Z^{\mathrm{+}}},\; m,\, n\:\geqq\:0$
with $Q$ being the computable \emph{nucleus} value.

\medskip{}

And a\emph{ }colide-15 we will refer\emph{ }to by equation \ref{eq:COLIDE-15}
below :

\[
Q\,=\mathbb{\mathit{\,}\mathit{\left(6\times m\times n\right)}}\:+\:\mathit{\mathit{\left(5\times m\right)}+\:\mathit{n\:\:}\:\;\quad}\qquad\qquad\ref{eq:COLIDE-15}
\]

where $\; m,\, n\:$ are the \textbf{HI-level}s of the two constituent
primes/composites with $\; m,\, n\:\in\:\mathbb{Z^{\mathrm{+}}},$
$m\,>\:0$ , $n\:\geqq\:0$ with $Q$ being the readily computable
\emph{nucleus\index{Nucleus}} value.

\medskip{}

The importance of these three colides is that \emph{\uline{they
have no integer solutions}} for primes, else if they possibly do,
by contradiction they are composite as we can derive back $m$ and
$n$ which will be the Hexile levels of occurence of its two multiplicand
integers.

Using these equations in the following section \ref{sec:Sequencing-of-Nuclei-Values-H1_H5 Composites},
we will attempt to examine the sequencing of \emph{nuclei} values
in these colides : which of course only references the Hexile level
of composites.

\rule[0.5ex]{1\columnwidth}{1pt}

\section{Sequencing of Nuclei values for \emph{$\mathbb{\mathit{\mathbb{H^{\mathrm{1}}}}}$
}and \emph{$\mathbb{\mathit{\mathbb{H^{\mathrm{5}}}}}$ }\label{sec:Sequencing-of-Nuclei-Values-H1_H5 Composites}\emph{.}}

In the following sub-sections, the sequencing of composites in \emph{$\mathbb{\mathit{\mathbb{H^{\mathrm{1}}}}}$
}and \emph{$\mathbb{\mathit{\mathbb{H^{\mathrm{5}}}}}$ }is reviewed\emph{
}for the obvious non-random nature of its series\emph{.}

\subsection{Evaluating composite occurrences for a Colide-11 composite.\label{sec:Evaluating-composite-occurrences-colide-11-equation}}

Recalling a  colide-11 viz. :

\[
\mathbb{\mathit{Q\,=\,}\mathit{\left(6\times m\times n\right)+\:\mathit{m+\:\mathit{n\:}\:\;\quad}}}\qquad\qquad\eqref{eq:COLIDE-11}
\]

where $\; m,\, n\:\in\:\mathbb{Z^{\mathrm{+}}}$ are the \textbf{HI-level}s
of the two constituent primes/composites, with $Q$ being the readily
computable \emph{nucleus} value.

\medskip{}

We can transmute equation (\ref{eq:COLIDE-11}) into the following
function : 
\begin{equation}
f_{11}\left(\, m,\, n\,\right)\,=\left[\left(6\times m\times n\right)\right]+\:\mathit{m+\:\mathit{n\:}\:\;\quad}\label{eq:Colide-11-general-Function}
\end{equation}

this function taking on unordered integer pairs : $\left(\, m,\, n\,\right)$
of the\textbf{ }HI-levels of the two constituent primes/composites,
where $\; m,\, n\:\in\:\mathbb{Z^{\mathrm{+}}}$, $\; m,\, n\:>\:0$.

We will designate this function as the \uline{c}olide-11\textbf{
}\uline{s}tate function, abbreviated to CS-11  function.\medskip{}

Table \ref{tab:Table_Matrix_H1_p1q1} shows the values obtained by
applying to the CS-11 function the first 8 positive integers excluding
zero. The resulting symmetric matrix, with its main diagonal is highlighted
in bold.

\textbf{\huge }
\begin{table}[H]
\textbf{\huge }%
\begin{tabular}{|c|cccccccc|l|}
\hline 
\multicolumn{10}{|c|}{}\tabularnewline
\hline 
\hline 
$\mathrm{LEVEL\,:\left(\mathrm{m}\right):\,\,\,\,\,}$ &  & 1 & 2 & 3 & 4 & 5 & 6 & 7 & 8\tabularnewline
\hline 
$\mathit{f_{11}\left(\, m,\, n\,\right)=\left[\left(6\times m\times n\right)\right]+m+\mathit{n}}$ &  & $\mathit{c_{N}^{11}\,}$  & $\mathit{c_{N}^{11}\,}$  & $\mathit{c_{N}^{11}\,}$  & $\mathit{c_{N}^{11}\,}$  & $\mathit{c_{N}^{11}\,}$  & $\mathit{c_{N}^{11}\,}$  & $\mathit{c_{N}^{11}\,}$  & $\mathit{c_{N}^{11}\,}$ \tabularnewline
\hline 
$\mathrm{LEVEL\,:\left(\mathrm{n}\right):\,\,\,\,\,}$ &  &  &  &  &  &  &  &  & \tabularnewline
$\mathrm{n\,=\,}$1 &  & \textbf{\emph{8}} & \emph{15} & \emph{22} & \emph{29} & \emph{36} & \emph{43} & \emph{50} & \emph{57}\tabularnewline
$\mathrm{n\,=\,}$2 &  & \emph{15} & \textbf{\emph{28}} & \emph{41} & \emph{54} & \emph{67} & \emph{80} & \emph{93} & \emph{106}\tabularnewline
$\mathrm{\mathrm{n\,=\,}}$3 &  & \emph{22} & \emph{41} & \textbf{\emph{60}} & \emph{79} & \emph{98} & \emph{117} & \emph{136} & \emph{155}\tabularnewline
$\mathrm{n\,=\,}$4 &  & \emph{29} & \emph{54} & \emph{79} & \textbf{\emph{104}} & \emph{129} & \emph{154} & \emph{179} & \emph{204}\tabularnewline
$\mathrm{n\,=\,}$5 &  & \emph{36} & \emph{67} & \emph{98} & \emph{129} & \textbf{\emph{160}} & \emph{191} & \emph{222} & \emph{253}\tabularnewline
$\mathrm{n\,=\,}$6 &  & \emph{43} & \emph{80} & \emph{117} & \emph{154} & \emph{191} & \textbf{\emph{228}} & \emph{265} & \emph{302}\tabularnewline
$\mathrm{n\,=\,}$7 &  & \emph{50} & \emph{93} & \emph{136} & \emph{179} & \emph{222} & \emph{265} & \textbf{\emph{308}} & \emph{351}\tabularnewline
$\mathrm{n\,=\,}$8 &  & \emph{57} & \emph{106} & \emph{155} & \emph{204} & \emph{253} & \emph{302} & \emph{351} & \textbf{\emph{400}}\tabularnewline
\hline 
... &  & ... & ... & ... & ... & ... & ... & ... & ...\tabularnewline
\hline 
\end{tabular}\textbf{\huge \caption{\label{tab:Table_Matrix_H1_p1q1}Table showing the resulting $\mathrm{8\,\times\,8}$
- symmetric matrix of\textbf{ }colide-11 values obtained on applying
the CS-11 function to the various combination-pairs of the first 8
positive integers excluding zero.}
}
\end{table}
{\huge \par}

\medskip{}

The resulting tabulation of values in Table \ref{tab:Table_Matrix_H1_p1q1}
is a symmetric matrix, showing the predictable sequencing of nucleus
values for a colide-11 composite.. This symmetry is due to the commutative-mutable-emplacement
of the $m\,$ and $n\,$ variables, within the CS-11 function thus
:
\[
f_{11}\left(\, m,\, n\,\right)\,=\left(6\times m\times n\right)+\:\mathit{m+\:\mathit{n\:\:\;\quad\:\:\;\quad\:}\:\;\quad}\qquad\quad\eqref{eq:Colide-11-general-Function}
\]
\begin{equation}
i.e.\:\; f_{11}\left(\, m,\, n\,\right)\,=\left[\left(6\times m\right)\,+\,1\right]\times n+\:\mathit{m\,=\,\left[\left(6\times n\right)\,+\,1\right]\times m+\:\mathit{n}}\label{eq::Specific-Sequenced-CS-11-Function}
\end{equation}

The above equation \ref{eq::Specific-Sequenced-CS-11-Function} of
a\textbf{ }CS-11 function, thus governs the sequencing of colide-11
Hexile levels - with an arithmetic progression, hence producing a
series whose sequence is predictable and non-random in nature. 

\medskip{}

\shadowbox{\begin{minipage}[t]{1\columnwidth}%
Of critical note, is that the main diagonal of the symmetric matrix
in Table \ref{tab:Table_Matrix_H1_p1q1} represents those instances
where $m\,$=$\, n$ of our CS-11 function. 

To which, applying values of $m\,$=$\, n$ to equation \ref{eq:Colide-11-general-Function}
will confirm as computationally defined by:
\begin{equation}
f_{11}\left(\, m,\, m\,\right)\,=\left(6\,\times\, m^{2}\right)+\:\mathit{\left(2\,\times\, m\right)}\label{eq:CS-11-Function-main-Diagonal}
\end{equation}
\end{minipage}}

\rule[0.5ex]{1\columnwidth}{1pt}

\subsection{Evaluating composite occurrences for a Colide-55 composite.\label{sub:Evaluating-composite-occurrences-for-a-Colide-55-composite}}

We met in the section \ref{sec:Sequencing-of-Nuclei-Values-H1_H5 Composites}
the equation for a colide-55 viz.: 
\[
\mathbb{\mathit{Q\,=\,}\mathit{\left(6\times m\times n\right)+\:\mathit{\left(5\times m\right)+\:\mathit{\left(5\times n\right)\,+\,4\:}\:\;\quad}}}\qquad\qquad\eqref{eq:COLIDE-55}
\]

where $\; m,\, n\:\in\:\mathbb{Z^{\mathrm{+}}}$ are the HI-levels
of the two constituent \index{Primes}primes/composites, with $Q$
being the readily computable \emph{nucleus} value.

\smallskip{}

Whereby equation (\ref{eq:COLIDE-55}) can be transmuted to the function
:
\begin{equation}
f_{55}\left(\, m,\, n\,\right)\,=\mathit{\left(6\times m\times n\right)+\:\mathit{\left(5\times m\right)+\:\mathit{\left(5\times n\right)\,+\,4\:}\:\;\quad}}\label{eq:Colide-55-General-Function}
\end{equation}

this function taking on unordered integer pairs : $\left(\, m,\, n\,\right)$
of the\textbf{ }HI-levels of the two constituent primes/composites,
where $\; m,\, n\:\in\:\mathbb{Z^{\mathrm{+}}}$, with $\; m,\, n\:\geqq\:0$.

We will designate this function as the \uline{c}olide-55 \textbf{\uline{s}}tate
function, abbreviated to CS-55 function.\medskip{}

Table \ref{tab:Table-showing-the-CS-55-Function-values}~ shows the
values obtained by applying to the CS-55 function the first 8 positive
integers excluding zero.The resulting symmetric matrix, with its main
diagonal is highlighted in bold.

\textbf{\huge }
\begin{table}[H]
\textbf{\huge }%
\begin{tabular}{|c|c|cccccc|l|}
\hline 
\multicolumn{9}{|c}{}\tabularnewline
\hline 
\hline 
{\small $\mathrm{LEVEL\,:\left(\mathrm{m}\right):\,\,\,\,\,}$} & {\small 0} &  & {\small 1} & {\small 2} & {\small 3} & {\small 4} & {\small 5} & {\small 6}\tabularnewline
\hline 
{\small $\mathit{f_{55}\left(m,n\right)=\mathit{\left(6\times m\times n\right)+\left(5\times m\right)+\left(5\times n\right)+4}}$} & {\small $\mathit{c_{N}^{55}\,}$ } &  & {\small $\mathit{c_{N}^{55}\,}$ } & {\small $\mathit{c_{N}^{55}\,}$ } & {\small $\mathit{c_{N}^{55}\,}$ } & {\small $\mathit{c_{N}^{55}\,}$ } & {\small $\mathit{c_{N}^{55}\,}$ } & {\small $\mathit{c_{N}^{55}\,}$ }\tabularnewline
\hline 
{\small $\mathrm{LEVEL\,:\left(\mathrm{n}\right):\,\,\,\,\,}$} &  &  &  &  &  &  &  & \tabularnewline
\hline 
{\small $\mathrm{n\,=\,}$0} & \textbf{\small 4} &  & {\small 9} & {\small 14} & {\small 19} & {\small 24} & {\small 29} & {\small 34}\tabularnewline
{\small $\mathrm{n\,=\,}$1} & {\small 9} &  & \textbf{\small 20} & \emph{\small 31} & \emph{\small 42} & {\small 53} & {\small 64} & {\small 75}\tabularnewline
{\small $\mathrm{n\,=\,}$2} & {\small 14} &  & \emph{\small 31} & \textbf{\emph{\small 28}} & \emph{\small 41} & \emph{\small 54} & \emph{\small 67} & \emph{\small 80}\tabularnewline
{\small $\mathrm{\mathrm{n\,=\,}}$3} & {\small 19} &  & \emph{\small 42} & \emph{\small 41} & \textbf{\emph{\small 60}} & \emph{\small 79} & \emph{\small 98} & \emph{\small 117}\tabularnewline
{\small $\mathrm{n\,=\,}$4} & {\small 24} &  & {\small 53} & \emph{\small 54} & \emph{\small 79} & \textbf{\emph{\small 104}} & \emph{\small 129} & \emph{\small 154}\tabularnewline
{\small $\mathrm{n\,=\,}$5} & {\small 29} &  & \emph{\small 64} & \emph{\small 67} & \emph{\small 98} & \emph{\small 129} & \textbf{\emph{\small 160}} & \emph{\small 191}\tabularnewline
{\small $\mathrm{n\,=\,}$6} & {\small 34} &  & {\small 75} & \emph{\small 80} & \emph{\small 117} & \emph{\small 154} & \emph{\small 191} & \textbf{\emph{\small 228}}\tabularnewline
\hline 
... &  &  & ... & ... & ... & ... & ... & ...\tabularnewline
\hline 
\end{tabular}\textbf{\huge \caption{\label{tab:Table-showing-the-CS-55-Function-values}Table showing
the resulting $\mathrm{7\,\times\,7}$ - symmetric matrix of\textbf{
colide-55\index{Colide-55}} values obtained on applying the \textbf{CS-55\index{CS-55}}
function to the various combination-pairs of the first 7 positive
integers including zero.}
}
\end{table}
{\huge \par}

The resulting tabulation of values in Table \ref{tab:Table-showing-the-CS-55-Function-values}
is a symmetric matrix, showing the predictable sequencing of nucleus
values for a colide- 55 composite. This symmetry is due to the commutative-mutable-emplacement
of the $m\,$ and $n\,$ variables, within the CS-55 function thus
:

\[
f_{55}\left(\, m,\, n\,\right)=\mathit{\left(6\times m\times n\right)+\left(5\times m\right)+\left(5\times n\right)+\,4\:\:\;\quad}\qquad\qquad\quad\eqref{eq:Colide-55-General-Function}
\]
\begin{equation}
i.e.\:\; f_{55}\left(\, m,\, n\,\right)=\left[\left(6\times m\right)\,+\,5\right]\times n+\:\mathit{\left(5\times m\right)\,+\,4=\,\left[\left(6\times n\right)\,+\,5\right]\times m+\:\mathit{\left(5\times n\right)}+\,4}\label{eq::Specific-Sequenced-CS-55-Function}
\end{equation}
The above equation \ref{eq::Specific-Sequenced-CS-55-Function} of
a CS-55 function, thus governs the sequencing of colide-55\textbf{
}Hexile levels - with an arithmetic progression, hence producing a
series whose sequence is predictable and non-random in nature.. 

\medskip{}

\shadowbox{\begin{minipage}[t]{1\columnwidth}%
Also of critical note is that the main diagonal of the symmetric matrix
in Table \ref{tab:Table-showing-the-CS-55-Function-values} represents
those instances where $m\,$=$\, n$ of our \textbf{CS-55} function.

To which, applying values of $m\,$=$\, n$ to equation \ref{eq:Colide-55-General-Function}
will confirm as computationally defined by equation \ref{eq:CS-55-Main-Diagonal}
below:
\begin{equation}
f_{55}\left(\, m,\, m\,\right)\,=\left(6\,\times\, m^{2}\right)+\:\mathit{\left(10\,\times\, m\right)\,+\,4}\label{eq:CS-55-Main-Diagonal}
\end{equation}
\end{minipage}}

\rule[0.5ex]{1\columnwidth}{1pt}

\subsection{Evaluating composite occurrences for a Colide-15 composite.\label{sub:Evaluating-composite-occurrences-for-a-Colide-15 composite}}

Looking at the equation for a colide-15 from section \ref{sec:Sequencing-of-Nuclei-Values-H1_H5 Composites}
to wit:
\[
\mathbb{\mathit{Q\,=\,}\mathit{\left(6\times m\times n\right)+\,\left(5\times m\right)\,+\, n\:\:\;\quad}}\qquad\qquad\eqref{eq:COLIDE-15}
\]

Whereby equation (\ref{eq:COLIDE-15}) can be transmuted to the function
:
\begin{equation}
f_{15}\left(\, m,\, n\,\right)\,=\mathit{\left(6\times m\times n\right)+\:\mathit{\left(5\times m\right)+\:\mathit{n\,\:}\:\;\quad}}\label{eq:Colide-15-State-Function}
\end{equation}

this function taking on ordered integer pairs : $\left(\, m,\, n\,\right)$
of the\textbf{ }HI-levels of the two constituent primes/composites,
where $\; m,\, n\:\in\:\mathbb{Z^{\mathrm{+}}}$, with $\; m,\, n\:\geqq\:0$.

We will designate this function as the \uline{c}olide-15\textbf{
}\uline{s}tate function, abbreviated to CS-15 function.

Also of critical note is that the input values are ordered pairs $\left(\, m,\, n\,\right)$,
and that they do not necessarily represent the HI-levels\textbf{ }of
primes/composites.

\medskip{}

Table \ref{tab:Table-showing-the-CS-15-Returned-Values}~ shows the
values obtained by applying to the CS-55 function the first 8 positive
integers excluding zero.

\textbf{\huge }
\begin{table}[H]
\textbf{\huge }%
\begin{tabular}{|c|c|ccccccc|l|}
\hline 
\multicolumn{10}{|c}{}\tabularnewline
\hline 
\hline 
$\mathrm{LEVEL\,:\left(\mathrm{m}\right):\,\,\,\,\,}$ & 1 &  & 2 & 3 & 4 & 5 & 6 & 7 & 8\tabularnewline
\hline 
$\mathit{f_{15}\left(m,n\right)=\mathit{\left(6\times m\times n\right)+\left(5\times m\right)+n}}$ & $\mathit{c_{N}^{15}\,}$  &  & $\mathit{c_{N}^{15}\,}$  & $\mathit{c_{N}^{15}\,}$  & $\mathit{c_{N}^{15}\,}$  & $\mathit{c_{N}^{15}\,}$  & $\mathit{c_{N}^{15}\,}$  & $\mathit{c_{N}^{15}\,}$  & $\mathit{c_{N}^{15}\,}$ \tabularnewline
\hline 
$\mathrm{LEVEL\,:\left(\mathrm{n}\right):\,\,\,\,\,}$ &  &  &  &  &  &  &  &  & \tabularnewline
\hline 
$\mathrm{n\,=\,}$1 & 12 &  & 23 & 34 & 45 & 56 & 67 & 78 & 89\tabularnewline
$\mathrm{n\,=\,}$2 & 19 &  & 36 & \emph{53} & \emph{70} & 87 & 104 & 121 & 138\tabularnewline
$\mathrm{n\,=\,}$3 & 26 &  & \emph{49} & \emph{72} & \emph{95} & \emph{118} & \emph{141} & \emph{164} & 187\tabularnewline
$\mathrm{\mathrm{n\,=\,}}$4 & 33 &  & \emph{62} & \emph{91} & \emph{120} & \emph{149} & \emph{178} & \emph{207} & \emph{236}\tabularnewline
$\mathrm{n\,=\,}$5 & 40 &  & 75 & \emph{110} & \emph{145} & \emph{180} & 215 & \emph{250} & 285\tabularnewline
$\mathrm{n\,=\,}$6 & 47 &  & \emph{88} & \emph{129} & \emph{170} & 211 & \emph{252} & \emph{293} & 334\tabularnewline
$\mathrm{n\,=\,}$7 & 54 &  & 101 & \emph{148} & \emph{195} & \emph{242} & \emph{289} & \emph{336} & 383\tabularnewline
$\mathrm{n\,=\,}$8 & 61 &  & 114 & \emph{167} & 220 & 273 & \emph{326} & \emph{379} & 432\tabularnewline
\hline 
... &  &  & ... & ... & ... & ... & ... & ... & ...\tabularnewline
\hline 
\end{tabular}\textbf{\huge \caption{\label{tab:Table-showing-the-CS-15-Returned-Values}The resultant
$\mathrm{8\,\times\,8}$ - matrix table, showing the\textbf{ colide-15\index{Colide-15}}
values obtained on applying the \textbf{CS-15} function to various
combination-pairs of the first 8 positive integers excluding zero.}
}
\end{table}
{\huge \par}

The resulting tabulation of values in Table \ref{tab:Table-showing-the-CS-55-Function-values}
shows the predictable sequencing of nucleus values for a colide-15
composite.

\[
f_{15}\left(\, m,\, n\,\right)\,=\mathit{\left(6\times m\times n\right)+\:\mathit{\left(5\times m\right)+\:\mathit{n\,\:}\:\;\quad}}\qquad\eqref{eq:Colide-15-State-Function}
\]
\begin{equation}
i.e.\:\; f_{11}\left(\, m,\, n\,\right)\,=\left[\left(6\times m\right)\,+\,1\right]\times n+\:\mathit{\left(5\times m\right)\,=\,\left[\left(6\times n\right)\,+\,5\right]\times m+\:\mathit{n}}\label{eq:Specific-Sequenced-CS-15-Function}
\end{equation}

From equation \ref{eq:Specific-Sequenced-CS-15-Function} it is obvious
that the CS-15 function sequences \emph{nuclei} values is an arithmetic
progression, hence producing a series whose sequence is predictable
and non-random in nature.

\rule[0.5ex]{1\columnwidth}{1pt}

\subsection{Conclusion to Sequencing of Nuclei values for \emph{$\mathbb{\mathit{\mathbb{H^{\mathrm{1}}}}}$
}and \emph{$\mathbb{\mathit{\mathbb{H^{\mathrm{5}}}}}$ }elements\emph{.\label{sub:Conclusion-to-Sequencing-Nuclei-values-H1-H5 }}}

In the preceding sub-sections (\ref{sec:Evaluating-composite-occurrences-colide-11-equation})
and (\ref{sub:Evaluating-composite-occurrences-for-a-Colide-55-composite}),
we were able to establish that colide-11 and colide-55 composites
which populate \emph{$\mathbb{\mathit{\mathbb{H^{\mathrm{1}}}}}$,
}occur in a non-random sequence. Then from the result in section (\ref{sec:Distinguishing-feature-of-Composites_Primes})
where we assayed \emph{$\mathbb{\mathit{\mathbb{H^{\mathrm{1}}}}}$
}to consist of only\emph{ }colide-11 and colide-55 composites and
primes, this means we can logically deduce : that primes which are
the only other elements of \emph{$\mathbb{\mathit{\mathbb{H^{\mathrm{1}}}}}$,
}occur in a non-random fashion in the vacant Hexile levels not occupied
by colide-11s and colide-55s.

With similiar reasoning from the conclusion drawn in sub-section (\ref{sub:Evaluating-composite-occurrences-for-a-Colide-15 composite})
we found that the colide-15 composites which populate \emph{$\mathbb{\mathit{\mathbb{H^{\mathrm{5}}}}}$
}also occur in a non-random sequence within this equivalence class.
Then from the result in section (\ref{sec:Distinguishing-feature-of-Composites_Primes})
where we assayed \emph{$\mathbb{\mathit{\mathbb{H^{\mathrm{5}}}}}$
}to consist of only\emph{ }colide-15 composites and primes, this also
means we can logically deduce : that primes which are the only other
elements found in \emph{$\mathbb{\mathit{\mathbb{H^{\mathrm{5}}}}}$,}
also\emph{ }occur in a non-random fashion in the vacant Hexile levels
not occupied by colide-15s.

\smallskip{}

The conclusion can then be drawn that we have crudely established
the rhythm and harmonics in the symphony of prime and composite integers.

\smallskip{}

As a preliminary to the next section \ref{sec:Counting-primes-in-H-H5},
where the cardinality of primes in \emph{$\mathbb{\mathit{\mathbb{H^{\mathrm{1}}}}}$
}and \emph{$\mathbb{\mathit{\mathbb{H^{\mathrm{5}}}}}$} is examined,
we examine the following prime sieve. We introduce this, as the over-arching
theoretical underpinnings of being able to sieve primes, also bears
directly on being able to count them.

We found a way to generate the sequence of integers representing the
\emph{\uline{Hexile levels}} of composites populating \emph{$\mathbb{\mathit{\mathbb{H^{\mathrm{1}}}}}$}
using the CS-11 and CS-55 functions. This means we can conceptualise
and designate this series of integers as a set symbolized by $\mathrm{S}$.
Then the set $\mathrm{T}$, which is the set representing the \emph{\uline{Hexile
levels}} of primes in \emph{$\mathbb{\mathit{\mathbb{H^{\mathrm{1}}}}}$
}can be naively found by \emph{:
\begin{equation}
\mathrm{T}=\mathbb{N\,\mathrm{-}\,\mathrm{S}}\label{eq:Set-of-Primes-H1}
\end{equation}
}

Similiarly for primes in \emph{$\mathbb{\mathit{\mathbb{H^{\mathrm{5}}}}}$
we can} generate the set of integers representing composite Hexile
levels using the CS-15 function and designate this set by $\mathrm{V}$.
Then the set $\mathrm{W}$ which is the set representing the Hexile
levels of primes in \emph{$\mathbb{\mathit{\mathbb{H^{\mathrm{5}}}}}$
}can also naively be found by \emph{:
\begin{equation}
\mathrm{W}=\mathbb{N\,\mathrm{-}\,\mathrm{V}}\label{eq:Set-of-Primes-H5}
\end{equation}
}

\rule[0.5ex]{1\columnwidth}{1pt}

\pagebreak{}

\section{Counting primes in \emph{$\mathbb{\mathit{\mathbb{H^{\mathrm{1}}}}}$
and $\mathbb{\mathit{\mathbb{H^{\mathrm{5}}}}}$.\label{sec:Counting-primes-in-H-H5}}}

With the three colides serving as the theoretical basis for composites,
we have gained a little insight into the multiplicative mechanics
of composite generation.%
\footnote{Not all composites in $\mathbb{N}$ are in our present consideration,
but composites which are important to asymmetric encryption in the
PKI industry: as found in \emph{$\mathbb{\mathit{\mathbb{H^{\mathrm{1}}}}}$
}and\emph{ $\mathbb{\mathit{\mathbb{H^{\mathrm{5}}}}}$}. For composites
in equivalence classes : \emph{$\mathbb{\mathit{\mathbb{H^{\mathrm{0}}}}},$
$\mathbb{\mathit{\mathbb{H^{\mathrm{2}}}}},$ $\mathbb{\mathit{\mathbb{H^{\mathrm{3}}}}}$}
and \emph{$\mathbb{\mathit{\mathbb{H^{\mathrm{4}}}}}$} were not explored,
as these HICs do not contain primes greater than 5, and less importantly
their composites are easily decomposed into multiples of either 2
or 3.%
}

Then it should be easy for us to find a process whereby and wherein
for \emph{any} given positive integer\emph{ }in\emph{ $\mathbb{\mathit{\mathbb{N}}}$,}
to undertake the following steps in computing the number of Hexile
levels occupied by primes less than that integer : 
\begin{enumerate}
\item Obtaining the nucleus value $Q$ for a given integer $c$.
\item Obtaining maximal Hexile level values which satisfies the three colide(s)
formulations viz:

\begin{enumerate}
\item The maximal Hexile level value of a colide-11, being the floor integer
which satisfies the colide-11, when one of its variables is equal
to 1, in the case of $m=1$ we get:
\[
\mathbb{\mathit{Q\,\geq\,}\mathit{\left(6\times m\times n\right)+\:\mathit{m+\:\mathit{n\:}\:\;\quad}}}\qquad\qquad
\]
\[
\mathbb{\mathit{\quad\,\geq\,}\mathit{\left(6\times1\times n\right)+\:\mathit{1+\:\mathit{n\:}\:\;\quad}}}\qquad\qquad
\]
\begin{equation}
\mathbf{\mathbb{\mathit{\therefore\, Q\,\geq\,}\mathit{\left(7\times n\right)+\:\mathit{1\;\qquad\qquad\;\:\:\;\quad}}}\qquad\qquad}\label{eq:Minimal-m-colide-11}
\end{equation}
Then $n_{11}\:$ the maximal Hexile Level value satisfying equation
(\ref{eq:Minimal-m-colide-11}) for a colide-11 is computed from:
\begin{equation}
n_{11}\,=\, floor\left[\frac{Q\,-\,1}{7}\right]\label{eq:n11-Maximal-colide-11}
\end{equation}
where {}``$\, floor\left[\:\;\right]$ '', refers to the integer
floor function, with $n_{11}\,\in\:\mathbb{Z^{\mathrm{+}}}.$
\item The maximal Hexile level value of a colide-55 being the floor integer
which satisfies the colide-55, when one of its variables is equal
to 1, in the case of $m=1$ we get:
\[
\mathbb{\mathit{Q\,\geq\,}\mathit{\left(6\times m\times n\right)+\:\mathit{\left(5\times m\right)+\:\mathit{\left(5\times n\right)+\,4\:}\:\;\quad}}}\qquad\qquad
\]
\[
\mathbb{\mathit{\quad\,\geq\,}\mathit{\left(6\times1\times n\right)+\:\mathit{5+\:\mathit{\left(5\times n\right)+\,4\:}\:\;\quad}}}\qquad\qquad
\]
 
\begin{equation}
\mathbf{\mathbb{\mathit{\therefore\, Q\,\geq\,}\mathit{\left[\left(11\times n\right)\right]+\:\mathit{9\qquad\qquad\qquad\qquad\qquad\qquad\;}}}\qquad\qquad}\label{eq:Minimal-Range-c-55}
\end{equation}
Then $n_{55}$, the maximal Hexile Level value satisfying equation
(\ref{eq:Minimal-Range-c-55}) for a colide-55 is computed from :
\begin{equation}
n_{55}\,=\, floor\left[\frac{Q\,-\,9}{11}\right]\label{eq:c-55-nup}
\end{equation}
where {}``$\, floor\left[\:\;\right]$ '', refers to the integer
floor function, with $n_{55}\,\in\:\mathbb{Z^{\mathrm{+}}}.$
\item The maximal Hexile level value of a colide-15, being the floor integer
which satisfies the colide-15 when one of its variables is equal to
1, in the case of $m=1$ we get :
\[
\mathbb{\mathit{Q\,\geq\,}\mathit{\left(6\times m\times n\right)+\:\mathit{\left(5\times m\right)+\:\mathit{n\,}}\:\:}}\qquad\qquad
\]

\end{enumerate}

\[
\mathbb{\mathit{\quad\,\geq\,}\mathit{\left[\left(6\times1\times n\right)\right]+\:\mathit{5+\: n\,\:\;\quad}}}\qquad\qquad
\]
\begin{equation}
\mathbb{\mathit{\therefore\, Q\,\geq\,}\mathit{\left[\left(7\times n\right)\right]+\:5\;\qquad\qquad\;\:\:\;\quad}}\qquad\label{eq:Minimal-range-c15-n5}
\end{equation}

Then $n_{15}\:$ the maximal Hexile Level value satisfying equation
equation equation \ref{eq:Minimal-range-c15-n5} for a colide-55 is
computed from :
\begin{equation}
n_{15}\,=\, floor\left[\frac{Q\,-\,5}{7}\right]\label{eq:c-15-nup}
\end{equation}
where {}``$\, floor\left[\:\;\right]$ '', refers to the integer
floor function, with $n_{15}\,\in\:\mathbb{Z^{\mathrm{+}}}$

\item The next stage is computing the number of composites below the given
integer $c$. It is here we can apply combinatorial or numerical analysis
for determining the \emph{number} of integer tuples $\left(\, m,\, n\,\right)$
satisfying the following equalities within the domain of integers
from 1 to their maximal Hexile Level(s) for the three 'flavours' of
colides. 

\begin{enumerate}
\item Where for the following colide-11 equation \ref{eq:Colide-11-Combinatorial-equation}
below : 
\begin{equation}
\mathbb{\mathit{Q\,\geq\,}\mathit{\left(6\times m\times n\right)+\:\mathit{m+\:\mathit{n\:}\:\;\quad}}}\qquad\qquad\label{eq:Colide-11-Combinatorial-equation}
\end{equation}
we compute the \emph{number} of integer tuples $\left(\, m,\, n\,\right)$
within the domain of integers in the set $\left\{ 1,...,n_{11}\right\} $
which satisfies the above equation \ref{eq:Colide-11-Combinatorial-equation}.
\\
Noting that we have to take into account instances where there are
\emph{\uline{duplicated Hexile level values}}. This resulting from
instances where the values of tuples $\left(m,\, n\right)$ are interchanged,
yielding the same Hexile level value in the CS-11 function. This occurs
through the commutatively emplaced variables as shown in equation
(\ref{eq::Specific-Sequenced-CS-11-Function}) below : 
\[
f_{11}\left(\, m,\, n\,\right)\,=\left[\left(6\times m\right)\,+\,1\right]\times n+\:\mathit{m\,=\,\left[\left(6\times n\right)\,+\,1\right]\times m+\:\mathit{n}}\qquad\qquad\qquad\qquad\qquad\eqref{eq::Specific-Sequenced-CS-11-Function}
\]
 Additionally we are to take into account \emph{\uline{none-duplicated}}
\emph{\uline{Hexile level values}}. This resulting from instances
where the values of tuples $\left(m,\, n\right)$ are identical. This
was noted before at the end of sub-section (\ref{sec:Evaluating-composite-occurrences-colide-11-equation})
as accounting for the symmetric matrix of graphed CS-11 values as
shown in equation (\ref{eq:CS-11-Function-main-Diagonal}) repeated
below : 
\[
f_{11}\left(\, m,\, m\,\right)\,=\left(6\,\times\, m^{2}\right)+\:\mathit{\left(2\,\times\, m\right)}\qquad\qquad\qquad\qquad\eqref{eq:CS-11-Function-main-Diagonal}
\]

\item Where for the following colide-55 equation \ref{eq:Colide-55-Combinatorial-Equation}
below : 
\begin{equation}
\mathbb{\mathit{Q\,\geq\,}\mathit{\left(6\times m\times n\right)+\:\mathit{\left(5\times m\right)+\:\mathit{\left(5\times n\right)+\,4}}}}\label{eq:Colide-55-Combinatorial-Equation}
\end{equation}
we compute the \emph{number} of integer tuples $\left(\, m,\, n\,\right)$
within the domain of integers in the set $\left\{ 1,...,n_{55}\right\} $
which satisfies the above equation \ref{eq:Colide-55-Combinatorial-Equation}.
\\
Noting that we have to take into account instances where there are
\emph{\uline{duplicated Hexile level values}}. This resulting from
instances where the values of tuples $\left(m,\, n\right)$ are interchanged,
yielding the same Hexile level value in the CS-55 function. This occurs
through the commutatively emplaced variables as shown in equation
(\ref{eq::Specific-Sequenced-CS-55-Function}) below : 
\[
\:\: f_{55}\left(m,n\right)=\left[\left(6\times m\right)+5\right]\times n+\left(5\times m\right)+4\qquad\qquad\qquad
\]
\[
\qquad\qquad\qquad=\left[\left(6\times n\right)+5\right]\times m+\left(5\times n\right)+4\qquad\quad\qquad\eqref{eq::Specific-Sequenced-CS-55-Function}
\]
Additionally we are to take into account \emph{\uline{none-duplicated}}
\emph{\uline{Hexile level values}}. This resulting from instances
where the values of tuples $\left(m,\, n\right)$ are identical. This
was noted before as accounting for the symmetric matrix of graphed
CS-55 values at the end of sub-section (\ref{sub:Evaluating-composite-occurrences-for-a-Colide-55-composite})
in equation (\ref{eq:CS-55-Main-Diagonal}) repeated below : 
\[
\:\qquad f_{55}\left(\, m,\, m\,\right)\,=\left(6\,\times\, m^{2}\right)+\:\mathit{\left(10\,\times\, m\right)\,+\,4}\qquad\qquad\qquad\eqref{eq:CS-55-Main-Diagonal}
\]
\\
and\\
\\

\item Where for a colide-15 equation \ref{eq:Colide-15-Combinatorial-Equation}
below : 
\begin{equation}
\mathbb{\mathit{Q\,\geq\,}\mathit{\left(6\times m\times n\right)+\:\mathit{\left(5\times m\right)+\:\mathit{n\,}}\:\:}}\qquad\qquad\label{eq:Colide-15-Combinatorial-Equation}
\end{equation}
we compute the \emph{number} of integer tuples $\left(\, m,\, n\,\right)$
within the domain of integers in the set $\left\{ 1,...,n_{15}\right\} $
which satisfies the above equation \ref{eq:Colide-15-Combinatorial-Equation}. 
\end{enumerate}
\item The number of primes can then be directly computed as the difference
of twice the nucleus value $Q$, less the total number of computed
composites as determined from the total number of 'viable' integer
tuples $\left(\, m,\, n\,\right)$ computed in the previous step.\\
We of course have to add into this figure the first three primes.
\end{enumerate}
In a general sense the essential logic in the above algorithm, is
predicated primarily on being able to precisely compute the total
number of $\mathbb{H^{\mathrm{1}}}$ and $\mathbb{H^{\mathrm{5}}}$
composites below any positive integer. This together with being able
to determine the total number of elements in $\mathbb{H^{\mathrm{1}}}$
and $\mathbb{H^{\mathrm{5}}}$ (as being twice the nucleus value $Q$
). Then we can determine by default and deduction both logically and
computationally with a high degree of assurance, the total number
of primes below that given integer. 

\rule[0.5ex]{1\columnwidth}{1pt}

\section{Conclusion to the Hexile Sieve Analysis of primes and composites.}

Even though our discussion is a basic summary of some of the areas
researched, three objectives were intended in this expression, namely
:
\begin{enumerate}
\item A simple explanation for critical mysteries surrounding primes and
composites.
\item A viable alternative observational and theoretical basis for re-examining
how Natural numbers are conceptualised.
\item And hopefully a stimulus and starting point to solving other key areas
of conjectures and interest in Number Theory.
\end{enumerate}
As for succesful intent and impact by these three aspirations into
prime number theory, and criticism into the overall and/or specific
logic and flaws in this discourse, the author welcomes and entertains
feedback.

\rule[0.5ex]{1\columnwidth}{1pt}

\rule[0.5ex]{1\columnwidth}{1pt}

\bigskip{}

\hspace{4cm}© Roger M. Creft, 2012.

\hspace{3.7cm}email : hexilehelix@yahoo.com
\end{document}